\documentclass[10pt,a4paper,reqno]{amsart}
\usepackage{amsfonts,amsthm,latexsym,amsmath,amssymb,amscd,amsmath, epsf}

\newcommand \la {\lambda}
\newcommand \pa {\partial}
\newcommand \bR {\mathbb R}
\newcommand \bC {\mathbb C}
\newcommand \CV {\mathcal V}
\newcommand \Si {\Sigma}
\newcommand \C {\mathbb C}
\newcommand \fL {\mathfrak {Lie}}
\newcommand \ka {\kappa}
\newcommand\R{{\mathbb R}}

\newcommand\CA{{\mathcal A}}

\newcommand\OO {\Omega^0}
\newcommand\one{\hbox{1\hskip -4pt 1}}

\begin{document}
             \numberwithin{equation}{section}

             \title[Giambelli-type  formula ]
             {Giambelli-type  formula for \\ subbundles of the tangent bundle}

	  \author[M.\`E.Kazarian]{Maxim Kazarian}
\address{Steklov Mathematical Institute, 42 Vavilova St., 117966, Moscow GSP-1, Russia} 
\email{kazarian@ium.ips.ras.ru}
\author[B.~Shapiro]{Boris Shapiro}
\address{Department of Mathematics, Stockholm University, SE-106 91, Stockholm, Sweden}
\email{shapiro@math.su.se}
\date{\today}

\begin{abstract}
Consider a generic $n$-dimensional subbundle $\CV$ of the tangent
bundle $TM$ on some given manifold $M$.
Given $\CV$ one can define different degeneracy loci $\Si_{r}(\CV),\;
r=(r_1\leq r_2\leq r_3
\leq\dots\leq r_k)$ on $M$ consisting of all points $x\in M$ for which the
dimension of the subspace $\CV^j(x)\subset TM(x)$ spanned
by all length $\leq j$ commutators of  vector fields tangent to $\CV$ at
$x$ is less than or equal to $r_j$.
 Under a certain  transversality assumption we 'explicitly' calculate  the $\mathbb Z_2$-cohomology classes of $M$ dual to $\Si_{r}(\CV)$ 
using determinantal formulas due to W.~Fulton and the expression for the Chern classes of the
associated bundle of the free Lie algebras in terms of the Chern classes of $\CV$.
\end{abstract}

\maketitle


\section {Preliminaries and results}
\label{sec:0}

\subsection {History and motivation} The question of the existence of a nontrivial  sub-bundle of the 
tangent bundle on a given manifold is a geometric problem of a long-standing interest. (Such sub-bundles are often called {\em distributions} and we will freely use both terms below.) In the basic nontrivial  case of rank  $2$ sub-bundles  first important results in the area go back to the classical treatise \cite {HH}.  Apparently the best achievements in this problem were obtained in late 60's by 
E.~Thomas in \cite {Th1},\cite {Th3}, see also a well-written survey \cite {Th2}.  Later some of his results were 
re-discovered by Y.~Matsushita, \cite{Ma}. Not much has been done in this area since then. One of the few recent exceptions is \cite {JM}. A rather detailed information is available about the existence of  (oriented)  sub-bundles of rank $2$. For rank $3$ and higher only the first obstruction to the latter problem is known, see \cite {Th3}. The algebraic invariants for this result  come from the Stiefel-Whitney classes of elements of $\mathbb{K}\tilde{\mathbb{O}}(M)$ which is the reduced real $K$-theory group of the manifold in question. Starting with late 70's the interest in the geometric properties of sub-bundles of the tangent bundle very stimulated by the development  of the singularity theory and the revival of the interest in nonholonomic mechanics.  A nice source of  information about this topic is \cite {Mont}.  In particular,  if a given manifold admits a sub-bundle of rank at least $2$ one can construct at each point  of the manifold an incomplete flag whose $i$th subspace is the linear span of the  commutators of length at most $i$  of the vector fields tangent to the sub-bundle. The ranks of these subspaces  will (in general) depend on the point, see below. For a small generic perturbation of the original sub-bundle the ranks of the subspaces of the incomplete flags will stay constant almost everywhere on the manifold and depend only on the rank of the original sub-bundle and the dimension of the manifold. Generalizing the question we started with one can formulate the following problem 

\medskip 
\noindent
\centerline {{\bf Problem.}{\em When does a manifold admit a distribution whose associated flags have }} 
{\em constant (and maximal possible) ranks throughout the manifold?} 

\medskip 
Being in general even more difficult than (still unsolved) initial question the latter problem has a nice answer in the case of oriented rank $2$-distributions on oriented $4$-manifolds. A rank $2$-distribution on a $4$-manifold whose associated flag has the set of ranks $(2,3,4)$ at each point is called the {\em Engel} distribution. Such a distribution has remarkable properties in many different aspects, see e.g. \cite {Ge1}-\cite{Ge3}. Namely, combining the results of \cite{Ge1} (see also \cite {KMS}) together with the recent \cite {Vo}, \cite {Vo1} we get

\medskip
\noindent
     {\bf Theorem.} {\sl An orientable $4$-manifold admits an orientable Engel distribution if and only if the manifold is parallellizable.}
     
     \medskip 
     The aim of the present paper is to develop the basics of the obstruction theory adjusted to the above  problem. 

\subsection{Standard definitions and notation}

The following  notions are standard, see e.g. \cite{GV}.
Let $M$ be an $m$-dimensional manifold and  $\CV\subset TM$ be an
$n$-dimensional sub-bundle (a rank $n$-distribution) 
on $TM$. 
 Given $\CV$ one associates at each point $x\in M$  its derived flag $f\CV_x=
\{\CV_x=\CV^1_x\subseteq \CV_x^2\subseteq\dots
\subseteq \CV^j_x\subseteq \dots\}, $ where
$\CV^j_x=\{\CV^{j-1}_x+[\CV^{j-1},\CV]\}_x$. If at each point
$x\in M$ there exists a positive integer $k(x)$ such that 
the  subspace $\CV^{k(x)}_x$ coincides with  $TM_x$  then $\CV$ is called
(maximally)
{\it nonholonomic.} Let
$n_j(x)$ denote the dimension of $\CV^j_x$. The set of numbers $(n_1(x),..., n_k(x))$ is called {\it the growth vector} of $\CV$ at $x$.
 For a given nonholonomic
$\CV$ the minimal number $k$ such that $\CV^k_x=T_xM$ at all points $x$
is called {\it the degree of nonholonomicity. }
(A sub-bundle $\CV$ is called {\it regular} if its
$n_j(x)$ do not depend on
$x$ and the corresponding set of numbers $n_1,n_2,\dots n_k$ is called
{\it the growth vector} of a regular sub-bundle $\CV$.)
Let $\mathfrak{Lie}_n$ denote the free Lie algebra with $n$ generators and let
$\mathfrak{Lie}_n^k$ be
its linear subspace  spanned by all elements of  length $k$.
 Let $d(n,k)$ be the dimension of
$\mathfrak{Lie}_n^k$.
  An  $n$-dimensional
$\CV$ with the degree of nonholonomicity $k$ is called {\it a maximal growth
sub-bundle} or {\it a mg-distribution}  if $n_j=\sum_{i\leq j}d(n,j)$ for all
$j<k$ and $n_k=m$. The growth vector with the entries $(n=d(n,1), d(n,1)+d(n,2), ..., d(n,1)+d(n,2)+...+d(n,k), m)$ is called the {\em maximal growth vector}.

\medskip 
\noindent
{\bf Remark. }  According to \cite {Ge2}
a germ of a generic distribution has maximal growth, where 'generic' means
belonging to some
open everywhere dense subset in $C^\infty$-Whitney topology. Thus locally a
typical sub-bundle is
a  mg-distribution while globally there are (many) topological obstructions to
 the existence of mg-distributions on a given $M$. The problem we are addressing in the present paper can be reformulated as constructing obstructions to the existence of mg-distributions on a given manifold. 

\medskip 
\noindent
{\bf Examples.} A contact structure is a regular mg-distribution.
A $2$-dimensional mg-subbundle on $M_4$ is called
{\it an Engel distribution}, see above and \cite {Ge1}.
It is the only example besides  contact structures and their
even-dimensional analogs with the
stable local normal form.

\subsection{Degeneracy loci}

Given a generic $n$-distribution $\CV\subset TM$ one expects that  globally
$M$ contains a (typically reducible) 
degeneracy locus $\Si$  consisting
of all such points $x$  where the growth vector $(n_1(x),...,n_k(x))$  is
 lexicographically smaller than the maximal one. Given a growth vector $ r=(r_1\leq r_2\leq
r_3\leq\dots\leq r_k)$ denote by  $\Si_{r}$ the subset of all $x\in M$ satisfying the
conditions $n_1(x)\leq r_1,...,
n_k(x)\leq r_k.$  Such
$\Si_{r}$ can be considered as  degeneracy loci in the standard
meaning of intersection theory, see \S 1, comp. \cite {F1}, ch. 14. Namely,
each $n$-distribution 
$\CV\subset TM$ induces the associated fiber bundle $\mathfrak{Lie}(\CV)\to M$
where the fiber
$\mathfrak{Lie}_x$ is the free Lie algebra generated by the subspace $\CV_x$.
$\mathfrak{Lie}(\CV)$
has an obvious grading $\mathfrak{Lie}(\CV)=\oplus_{k=1}^\infty
\mathfrak{Lie}^k(\CV)$ coming
from the grading $\mathfrak{Lie}_n=\oplus_{k=1}^\infty \mathfrak{Lie}_n^k$.
Moreover one can define a natural map  $\Phi: \mathfrak{Lie}(\CV)
\to TM$ of vector bundles sending each $\oplus_{i=1}^k\mathfrak{Lie}^i(\CV)$ onto $\CV^k$.
(The map $\Phi$ is not unique in very much the same way as the identifying map between $TM$ and its
nilpotentization, comp. \cite {GV}.)  This allows us to apply to the map $\Phi$
the determinantal formulas of \cite {F2}  under
the assumption that the considered $\Si_{r}$  has the expected
 (co)dimension, i.e. the same codimension as the corresponding
degeneracy locus  for a generic
map of flag bundles of the dimensions prescribed by $\mathfrak{Lie}_n$.  Algebraically, in order to be able to apply these formulas one also needs to express the  Chern classes of $\mathfrak{Lie}(\CV)$ in terms of  that of $\CV$.

\medskip 
\noindent 
{\bf Remark.} W.~Fulton in \cite {F2}Ê has generalized a large number of previously known special cases of determinantal formulas giving  the cohomology classes of different degeneracy loci for the maps of vector bundles to a very general situation of maps of flagged vector bundles. Such formulas could be traced back (through the works of Porteus and Thom as well as Laksov-Kempf and many other authors) to the pioneering results of G.~Giambelli on the degrees of  different strata in the spaces of matrices. For this reason in a number of publications the authors name similar determinantal formulas after Giambelli, see e.g. \cite {F2}, \S 7. For detailed account on determinantal formulas and degeneracy lovi we recommend \cite {FP} and for the information on Giovanni Giambelli see \cite {La}.

The contents of the paper is as follows. In $\S 2$ we
construct the map $\Phi: \mathfrak{Lie}(\CV)
\to TM$. In \S 3 we
find the explicit formula for the Chern
character of the bundle $\mathfrak {Lie}^k(\CV)$ in terms of the Chern
character of $\CV$ which turns out to be similar
to the formula for the dimensions of $\mathfrak {Lie}^k_n$. In principle, this
allows us in principle  to calculate the Chern classes of $\mathfrak
{Lie}^k(\CV)$ for any reasonable
specific example by inverting Newton polynomials, see Appendix. In $\S 4$ we recall the appropriate
determinantal  formula for maps of flag bundles, adjust it to our needs and
calculate some examples. Section \S 5 is devoted to some generalities
on derived flags $f\CV$ and the standard stratification of the spaces of matrices 
as well as counterexamples to transversality in big codimensions.
In $\S 6$ we enumerate all potentially admissible defect vectors
occurring for a generic distribution and prove the necessary transversality
 result showing that $\Si_{r}$ have the
expected (co)dimension in the cases $n=2,\; m\le 8$ and $n\ge 3, m\le 
\frac {n(n+1)(2n+1)}{6}$.  
In \S 7 we briefly  discuss some further directions of study and possible
generalizations of the transversality theorem.  The main result of the 
paper is formula $(11)$  justified  by the transversality theorem for the above mentioned 
values of $(n,m)$. Finally, appendix 
contains the  Mathematica code which explicitly calculates the 
necessary Chern classes of the homogeneous components of the free Lie 
algebra $\mathfrak{Lie}^k_{n}$ up to order $4$.  

In short, the main proposition of the paper can be summarized as 
follows.   In order to calculate the
universal formula for the cohomology class dual to $\Si_{r}$ for a given growth vector $r$ one has to
substitute the expressions for the Chern (or rather for the Stiefel-Whitney) 
classes of $\oplus_{i=1}^k\mathfrak {Lie}^i(\CV)$ and $TM$ into the appropriate
Giambelli-type formula.

The starting point of this note was the special case of $2$-sub-bundles on
$M_4$ considered in \cite {KMS}. Later the authors realized that analogous
computations over $\mathbb Z_2$ can be
carried out in the general setup described in the present paper. Sincere
thanks are due to R.~Montgomery  for important discussions of the subject and numerous remarks 
which substantially improved the quality of exposition. 
The second author is grateful to IHES(Paris), MPIM(Bonn) and 
MSRI(Berkeley) where parts of this project were carried out. 
\bigskip\bigskip

\section {$\Si_{r}$ as classical degeneracy loci} \label{sec:1}

\subsection {Associated bundle of free Lie algebras}

Given $\CV\subset TM$ let us now define the map $\Phi:\mathfrak {Lie}(\CV)\to
TM$ such that each sub-bundle
$\oplus_{i=1}^j\mathfrak {Lie}^i(\CV)$ is mapped onto $\CV^j$. In fact we
define another filtered
vector bundle $N(\CV)$ whose associated graded bundle  is isomorphic
to $\mathfrak {Lie}(\CV)$ with the canonical map $\Psi:N(\CV)\to TM$.
 As a result  one gets a (non-unique) map $\Phi:\mathfrak {Lie}(\CV)\to TM$ defined up to
a filtered isomorphism between  $N(\CV)$ and ${\mathfrak {Lie}}(\CV)$.

\subsection {Universal map}

     Let $\CV\subset TM$ be an $n$-sub-bundle  and $k$ be its degree
of nonholonomicity.

\medskip
\noindent
     {\bf Theorem 1.} {\sl There exists a  globally defined
 map of vector bundles
     $$\Phi_k:\oplus_{i\leq k}\fL^i(\CV)\to TM$$
     such that for all $x\in M$ and $j\le  k$,  the  subspace  $\CV^j(x)$
coincides with   the  image  of  $\Phi_j(x):\oplus_{i\leq
j}\fL^i(\CV)(x)\to  T_xM$,  where
$\Phi_j$ is the restriction of $\Phi_k$ to $\oplus_{i\leq
j}\fL^i(\CV)\subset \oplus_{i\leq k}\fL^i(\CV)$.}

     For the local version of this theorem, see \cite {GV}.

     {\sl Proof.} We construct an auxiliary flag of vector bundles
     $$N^1(\CV)\subset N^2(\CV)\subset\dots N^k(\CV)\,,$$
     the map $\Psi_k:N^k(\CV)\to TM$ satisfying the statement of Theorem,
and a canonical isomorphism
     $$N^j(\CV)/N^{j-1}(\CV)\cong \fL^j(\CV)$$
which shows that the flag of bundles $N(\CV)$ is isomorphic to the flag of
bundles $\fL(\CV)$.

Recall the notion of a Hall basis in the free Lie algebra $\mathfrak{Lie}_n$
with $n$ generators, see \cite {Bou}.
Namely, $\mathfrak{Lie}_n$  has  the following standard graded basis $H\,$ called
{\it Hall family} or {\it Hall basis}.
 Given a linearly ordered set $V$ (of cardinality $n$) let us define the
following
linearly ordered subset $H$  in the
free monoid  $Mon_V$.

1) if $u,v\in H$ and $lng(u)<lng(v)$ then $u<v$ where $lng$ denotes the usual
length of a word in $Mon_V$;
2) $V=H^1\subset H$ and $H^2$ consists of the set of all ordered pairs
$(v_1,v_2)\in H$ where $v_1<v_2$;
3) each element of $H$ of length at least 3 has the form
$a(bc)$ where $a,b,c\in H$, $bc\in H$, $b\leq a< bc$ and $b<c$. (Obviously,
$H=\bigcup_{k=1}^\infty
 H^k$ where $H^k$ is the set of all length $k$ elements in $H$.)

\noindent
{\bf Example.} If $V=\{u<v\}$ then $H^1=\{u;v\};$
$H^2=\{(u, v)\};$
$H^3=\{(u(u, v));$
$(v(u,v))\};$
$H^4=\{(u(u(u,v)));$
$(v(u(u,v)));$
$(v(v(u,v)))\};$
$H^5=$ $\{(u(u(u(u,v))));$
$(v(u(u(u,v))));$
$ (v(v(u$ $(u,v))));$
$ (v(v(v(u,v))));$
$ ((u,v)(u(u,v)));$
$((u,v)(v(u,v)))\}.$

  Now the construction  of  the flag of bundles $N(\CV)$ is as follows.
Let $W$ be the sheaf of free Lie algebras associated to the  sheaf  of
local  sections of the bundle $\CV\subset TM$.  The elements of $W$ are
$\R$-linear
combinations of Lie monomials of sections of $\CV$.  Denote  by  $\CA$  the
sheaf  of  rings of  smooth  functions  on  $M$.  Define  the  homomorphism
$D:W\to{\text der}\CA$ as follows.  If $v\in \mathcal V\subset W$ is of degree
1  we put
     $$D_v f=v \,f\,,\eqno(1)$$
  i.e.   the usual Lie derivative of the function $f$ along the vector
field $v$. Then, we assign by induction
     $$D_{[a,b]}f=D_b(D_af)-D_a(D_bf)\,.\eqno(2)$$

     The operation  $D$  is  well defined and $D_u$ is the derivation of
$\CA$ for every $u\in W$. Consider the sheaf $\CA\otimes_\R W$ of $\CA$-modules.
We introduce the Lie algebra structure on $\CA\otimes W$ as
     $$[f\otimes u, g \otimes w]=fg\otimes[u,w]+
                    fD_ug\otimes w-gD_wf\otimes u\eqno(3)$$
(the Jacobi identity is verified by direct computation). Finally, define
$N(\CV)$  to  be  the  quotient  Lie algebra of $\CA\otimes W$
by  the
ideal generated by all relations of the form
     $$f\otimes v=1\otimes fv\,,\qquad
                   f\in \CA\,,\quad v\in \CV\subset W\,.\eqno(4)$$

     Having (4) in mind  we drop the sign of the tensor product in  the
notation of the elements in  $N(\CV)$.  The filtration on $W$
 by the length of
Lie monomials gives the natural filtration $N^j(\CV)$  on
 the  sheaf   $N(\CV)$ of
$\CA$-modules.  We  claim that {\sl all $N^j(\CV)$ are locally free
sheafs of $\CA$-modules of finite ranks.} Indeed,  let $e_1,\dots,e_n$ be
the  set  of  local  sections of $\CV$ over some open domain $U\subset M$ such
that these sections form a basis in each fiber $\CV(x)$, $x\in U$. Then it
follows  from  (1)--(4)  that  every section $u$ of $W$ over $U$ can be
represented as
     $$u=\sum f_l h_l(e_1,...,e_n)\,,$$
     where $f_l$ are some functions and $h_l\in H$ are the elements of Hall
basis  of the free Lie algebra $\fL_n$.  Moreover,  this representation
is unique,  i.e. the set of sections $h_l(e_1,...,e_n)$,
$l\le\dim\oplus_{i\leq j}\fL^i_n$ forms the set of free generators of the
$\CA$-module $N^j(\CV)$. Thus the $\CA$-module $N^j(\CV)$ is  the  module
of  sections  of  some
vector bundle which we also denote  as $N^j(\CV)$.

     Observe, that   if   $[u,w]\in   N(\CV)$   has   degree   $j$   then
$[fu,gw]-fg[u,w]$ has degree strictly  less  than  $j$.  Therefore,  the
homomorphism $N^j(\CV)/N^{j-1}(\CV)\to \fL^j(\CV)$  is
well defined.  Moreover, the arguments above show that this homomorphism
is, in fact, an isomorphism of vector bundles.

     The homomorphism  of Lie algebras $\Psi:N(\CV)\to Vect(M)$
is now obvious.
It sends a  formal  Lie  bracket  of  vector  fields  in  $\CV$  to  the
corresponding commutator of these vector fields.  Formulas (1)--(4) show
that this    homomorphism    is    well    defined  and
$\Psi(N^j(\CV)(x))$ coincides with $\CV^j(x)$ by definition.
Theorem 1 is proved. \qed

\medskip 
\noindent 
    {\bf Remark.}  The  vector  bundle $N^j(\CV)$ can be also described,
as a usual
vector bundle by  trivializations and transition  functions.
Trivializations  of  $N^j(\CV)$ correspond to the trivializations of $\CV$
and are given by the sections $h_l(e_1,...,e_n)$,  $l\le\dim \oplus_{i\leq
j}\fL^i_n$,
where  $e_1,...,e_n$  are sections giving some local basis of $\CV$.  If
$\{e'_1,...,e'_n\}$ is
another basis such that $e_i=\sum a_{ir}e'_r$ then  to  find  transition
functions     for     $N^j(\CV)$    one    should    express    $h_l(\sum
a_{1p}e'_p,\dots,\sum a_{np}e'_p)$ using (1)--(4) as a linear combination
of $h_i(e'_1,...,e'_n)$ with some functional coefficients.

\bigskip\bigskip

\section {On the Chern classes of the bundle of free Lie
algebras} \label{sec:2}

     Let $E\to M$ be a complex vector bundle of dimension $n$ over a smooth
compact manifold $M$ (not necessarily a sub-bundle of $TM$). For any linear
representation  of the
group $GL(n,\bC)$ in $\C^m$  one can associate  in a natural way to the bundle $E\to M$ the
corresponding $m$-dimensional
bundle over $M$. For example, the bundles $E\otimes E$,
$E^*$, $\Lambda^2E$ etc.
are associated with the obvious representations of $GL(n,\bC)$ in $\C^n\otimes\C^n$,
$\C^{n*}$,
$\Lambda^2\C^n$ respectively. 

Given a basis  $\{e_1,\dots,e_n\}$ in $\C^n$ let  $\fL_n$ denote the free Lie algebra with the generators
$\{e_1,\dots,e_n\}$, and let  $\fL^k_n$ be its $k$th homogeneous component.
 A linear change of the above basis  acts naturally 
 on the spaces $\fL^k_n$. Denote by  $\fL^k_n(E)$ or simply   by $\fL^k(E)$
the bundle
over $M$ associated with this action.
The relation between the characteristic classes of the bundles $E$
and $\fL^k(E)$ is described in the following theorem. (This question was
already proposed
in \cite {Thr}.)

     Let $ch(E)\in H^*(M)$ be the Chern character of a bundle $E$. For any
element $\eta=\eta_0+\eta_1+\eta_2+\dots\in H^*(M)$ and a number $d$  set 
$(\eta)_d=\eta_0+\eta_1d+\eta_2 d^2+\dots$, where $\eta_i$ is homogeneous
of degree $2i$.

  \medskip 
  \noindent 
     {\bf Theorem 2}, see e.g. \cite {Re}. {\sl The Chern character of
the bundle
$\fL^k_n$ is given by
the  formula
  $$ch(\fL^k_n)={1\over k}\sum_{d|k}\mu(d)\Bigl(ch(E)^{k\over d}\Bigr)_d\;.\eqno(5)$$}

    Here the summation is taken over the set of all divisors of $k$ and
$\mu$ is the M\"obius
function. By taking the component of degree 0 in this formula we get the
well-known
expression for the dimension of $\fL^k_n$ (see [Bou])
  $$dim(\fL^k_n)={1\over k}\sum_{d|k}\mu(d)n^{k\over d}\;.$$

     {\sl Proof.} The main observation in this argument is that the total tensor algebra of
$\C^n$ is isomorphic,
as the $GL(n,\C)$-module, to the universal enveloping algebra of $\fL_n$.
Therefore, by
Poincar\'e-Birkhoff-Witt theorem one has
  $$T^*(E)\cong S^*(\fL(E))=S^*(\fL^1)\otimes S^*(\fL^2)\otimes\dots\;.$$

    Applying the Chern character to both sides of the latter equality we get
$${1\over1-ch(E)t}=\prod_{n=1}^\infty s(\fL^k)(t^k)\,,$$
where $t$ is a formal parameter, $s_i(V)=ch(S^iV)$ is  the Chern character of
the $i$th
symmetric power of $V$, and $s(V)$ is a formal series $s(V)(t)=\sum s_i(V)
t^i$. Now, applying
$-t\,d\log$ to both sides we get
  $${ch(E)t\over1-ch(E)t}=-\sum_{n=1}^\infty k\,t^k(d\log
s(\fL^k))(t^k)\,.\eqno(6)$$

Observe now that $-d\log s(V)(t)=(ch(V))_1+(ch(V))_2t+(ch(V))_3t^2+\dots$.
To prove
that we can assume (using the splitting principle) that the bundle $V$ is represented
as $V=V_1\oplus\dots\oplus V_m$, where the 1-dimensional bundle $V_i$ has
the Chern character $h_i$. Then
 $$-d\log s(V)=-d\log(s(V_1)s(V_2)\dots s(V_m))=-d\log\prod{1\over1-h_it}=
    \sum{h_i\over1-h_it}=$$ $$\qquad=
        \sum h_i+\sum h_i^2t+\sum
h_i^3t^2+\dots=(ch(V))_1+(ch(V))_2t+(ch(V))_3t^2\dots\,,$$
since $h_i^k=e^{kt_i}=\sum{t_i^lk^l\over l!}$, where $t_i$ denotes the first
Chern class of $V_i$.
Therefore, (6) is equivalent to
  $${ch(E)t\over1-ch(E)t}=
\sum(ch(\fL^1))_kt^k+2\sum(ch(\fL^2))_kt^{2k}+3\sum(ch(\fL^3))_kt^{3k}+\dots
\,,$$

       Comparing the  terms of the same degree in $t$ we get
  $$ch(E)^k=\sum_{d|k}d\,(ch(\fL^d))_{k\over d}\,.$$

        If we now multiply the $l$th homogeneous component of this equality
by $k^{-l}$ then
after this re-scaling we get
   $$(ch(E)^k)_{1\over k}=\sum_{d|k}d\,(ch(\fL^d))_{1\over d}\,.$$
  Applying to the latter equality  the M\"obius inversion formula we obtain
   $$k(ch(\fL^k))_{1\over k}=\sum_{d|k}\mu(d)(ch(E)^{k\over d})_{d\over k}\,.$$
which (after another re-scaling) gives the required formula of Theorem 2. \qed 

\medskip 
\noindent 
{\bf  Examples.} The relation between the Chern classes and the Chern character
gives the possibility to compute the Chern classes of $\fL^k_n$. For $k\le4$ (taking in account only terms of degree at most 4 in the 
characteristic classes) we obtain the following explicit formulas for the
total Chern class of $\fL^k_n$
 $$\vcenter{\halign{$#$\hfil\cr
c(\fL^1_n)= c(E)=1+c_1+c_2+c_3+c_4+\dots\,,\cr
c(\fL^2_n)= 1+( -1 + n ) \,c_1+
       \Bigl(( 1 - {{3\,n}\over 2} + {{{n^2}}\over 2} ) \,{{c_1}^2} + ( -2
+ n ) \,c_2\Bigr)+\cr
\qquad+\Bigl(( -1 + {{11\,n}\over 6} - {n^2} + {{{n^3}}\over 6} ) \,
       {{c_1}^3} + ( 4 - 4\,n + {n^2} ) \,c_1\,c_2 + ( -4 + n ) \,c_3\Bigr)+\cr
\qquad+\Bigl(
       (1 - {{25\,n}\over {12}} + {{35\,{n^2}}\over {24}} -
       {{5\,{n^3}}\over {12}} + {{{n^4}}\over {24}} ) \,{{c_1}^4} +
       ( -6 + 8\,n - {{7\,{n^2}}\over 2} + {{{n^3}}\over 2} )
\,{{c_1}^2}\,c_2 +\cr
\qquad\qquad+ ( 3 - {{5\,n}\over 2} + {{{n^2}}\over 2} )
       \,{{c_2}^2} + ( 9 - 6\,n + {n^2} ) \,c_1\,c_3 +
       ( -8 + n ) \,c_4\Bigr)+\dots\,,\cr
c(\fL^3_n)= 1+( -1 + {n^2} ) \,c_1+\Bigl(
       ( 2 - n - {{3\,{n^2}}\over 2} + {{{n^4}}\over 2} ) \,
       {{c_1}^2} + ( -3 + {n^2} ) \,c_2\Bigr)+\cr
\qquad+\Bigl(( -4 + 3\,n + {{17\,{n^2}}\over 6} -
       {n^3} - {n^4} +  {{{n^6}}\over 6} ) \,{{c_1}^3} +
       ( 12 - 4\,n - 5\,{n^2} + {n^4} ) \,c_1\,c_2 +
       ( -9 + {n^2} ) \,c_3\Bigr)+\cr
\qquad+\Bigl(( 8 - {{15\,n}\over 2} - {{61\,{n^2}}\over {12}} +
       {{7\,{n^3}}\over 2} + {{47\,{n^4}}\over {24}} - {{{n^5}}\over 2} -
       {{5\,{n^6}}\over {12}} + {{{n^8}}\over {24}} ) \,{{c_1}^4} +\cr
\qquad\qquad+( -36 + 19\,n + {{35\,{n^2}}\over 2} - 5\,{n^3} - 4\,{n^4} +
       {{{n^6}}\over 2} ) \,{{c_1}^2}\,c_2 + \cr
\qquad\qquad+( 18 - 6\,n - {{7\,{n^2}}\over 2} + {{{n^4}}\over 2} ) \,
       {{c_2}^2} + ( 36 - 6\,n - 11\,{n^2} + {n^4} ) \,c_1\,c_3 +
       ( -27 + {n^2} ) \,c_4\Bigr)+\dots\,,\cr
c(\fL^4_n)= 1+( -n + {n^3} )\,c_1+\Bigl(( 1 + n - {n^2} -
       {{{n^3}}\over 2} - {n^4} + {{{n^6}}\over 2} )
       \,{{c_1}^2} + ( -2\,(n + {n^3} )\,c_2\Bigr)+\cr
\qquad+\Bigl(( -4 - {n\over 3} + 2\,{n^2} + {{8\,{n^3}}\over 3} +
       {{3\,{n^4}}\over 2} - {n^5} - {{{n^6}}\over 2} - {{{n^7}}\over 2} +
       {{{n^9}}\over 6} )\,{{c_1}^3} + \cr
\qquad\qquad+( 8 + 4\,n - 4\,{n^2} - {n^3} - 3\,{n^4} + {n^6} )\,c_1\,
       c_2 + ( -4\,n + {n^3} )\,c_3\Bigr)+\cr
\qquad+\Bigl(( 13 - 2\,n - {{77\,{n^2}}\over {12}} - {{35\,{n^3}}\over 4} -
       {{3\,{n^4}}\over 4} + {{5\,{n^5}}\over 2} + {{55\,{n^6}}\over {24}} +
       {n^7} - {{{n^8}}\over 2} - {{{n^9}}\over 4} - {{{n^{10}}}\over 6} +
       {{{n^{12}}}\over {24}} )\,{{c_1}^4} + \cr
\qquad\qquad+( -48 + 12\,{n^2} + 18\,{n^3} + 7\,{n^4} - 5\,{n^5} -
       {{3\,{n^6}}\over 2} - 2\,{n^7} + {{{n^9}}\over 2} )\,
       {{c_1}^2}\,c_2 + \cr
\qquad\qquad+ ( 24 + 4\,n - 7\,{n^2} - {{{n^3}}\over 2} -
       2\,{n^4} + {{{n^6}}\over 2} )\,{{c_2}^2} +  \cr
\qquad\qquad+( 24 + 8\,n - 5\,{n^2} - {n^3} - 5\,{n^4} + {n^6} )\,c_1\,
       c_3 + ( -8\,n + {n^3} )\,c_4\Bigr)+\dots\,,\cr
}}$$

\medskip 
\noindent
{\bf Remark.} Substituting $w_i$ instead of $c_i$ in the above formulas and
reducing  coefficients mod 2 one gets the  expression for the
total Stiefel-Whitney
class of $\fL^k_n$ in the case of a real $n$-dimensional bundle $E$.
Note, that the coefficients of the above polynomials have integer values
for any $n$ and therefore their values mod $2$ are well defined.

\section {Determinantal formula and its application}\label{sec:3}

\subsection{Determinantal formula}
 First we recall a certain formula borrowed from \cite {F2}.
Assume that we have a flag $A_1\subset A_2\subset ...\subset A_l$ of the 
complex vector bundles over a manifold $M$ with the ranks $a_1\le a_2\le ... \le a_l$ resp. and a map
$$h:A_1\subset A_2\subset ...\subset A_l\to B$$
 to a manifold $B$ of dimension $b$. Assume furthermore that the set of  nonnegative integers
$\ka_1,...,\ka_l$ satisfies the inequalities 
$$0<a_1-\ka_1<a_2-\ka_2<...<a_l-\ka_l,\quad \ka_1<\ka_2<...<\ka_l<b.\eqno(7)$$

Let $\Omega_{\mathbb\ka}\subset M$ be the degeneracy locus defined by the conditions
$\text{rk} (h:A_i\to B)\le \ka_i,\;i=1,...,l.$ that is the set of all points $x\in M$ where all the previous conditions are valid. Now consider the Young diagram
$(p_1^{m_1},...,p_l^{m_l})$ where
$$p_1=a_l-\ka_l,\,p_2=a_{l-1}-\ka_{l-1},...,\,p_l=a_1-\ka_1,$$
$$m_1=b-\ka_l,\, m_2=\ka_l-\ka_{l-1},...,\,m_l=\ka_2-\ka_1.$$

Its dual diagram is $\mu=(q_1^{n_1},...,q_l^{n_l})$ where
$$q_1=b-\ka_1,\,...,\,q_l=b-\ka_l,$$
$$n_1=a_1-\ka_1,\,...,\,n_l=(a_l-\ka_l)-(a_{l-1}-\ka_{l-1}).$$

Let $cd(\mathbb \ka)=\vert \la \vert$, $s\la=b-\ka_1$, $s\mu=a_l-\ka_l$.
Finally, set
$$\rho(i)=max\{s\in [1,l]: i\le b-\ka_s=m_1+...+m_{l+1-s}\},\,i=1,...,s\la,$$
$$\rho^\prime(i)=min \{s\in [1,l]: i\le
a_s-\ka_1=n_1+...+n_{s}\},\,i=1,...,s\mu .$$

\medskip 
\noindent 
{\bf  Proposition 1}, see 10.2. of \cite {F2}. {\sl If the codimension of
$\Omega_{r}$
equals $cd(r)$ then the $\Bbb Z$-cohomology class  $[\Omega_{r}]_{\Bbb Z}$ of $M$ dual to $\Omega_{r}$ is given by
$$[\Omega_{r}]_{\Bbb Z}=\det(c_{\la_{i}-i+j}
(A_{\rho(i)}^*-B^*))_{1\le i,j\le s\la}=
\det(c_{\mu_{i}-i+j}(B-A_{\rho^\prime(i)}))_{1\le i,j\le s\mu},\eqno(8)$$
where $'*'$ denotes the dual bundle}.

\subsection {Real case} Consider a fixed flag of vector spaces
$$\bR^{a_1}\subset \bR^{a_2}\subset \dots \subset \bR^{a_l}.$$

Denote by $Mat(a_l,b)$ the space of all linear maps $\bR^{a_l}\to \bR^b$
identified with the space of $(a_l\times b)$-matrices. Two elements
$u,v\in Mat(a_l,b)$ are called {\it equivalent} if for all $i=1,...,l$ the
restriction of $u$ and $v$ to $\bR^{a_i}$ have the same rank. The set of all pairwise equivalent elements in $Mat(a_I,b)$ will be called a {\em stratum}. Obviously, one obtains  in this way  a finite stratification of  $Mat(a_I,b)$. 

Using the same notation as above consider now a map $h$ of {\it real} vector bundles

$$h:A_1\subset A_2\subset \dots \subset A_l\to B.$$

Let $x\in M$ be a point of the base and $U\ni x$ be its small neighborhood
such that the bundles are trivial and trivialized over $U$. Then the map $h$ over $U$ is
given by a family of matrices $$h_U:U\to Mat(a_l,b).$$

\noindent
{\bf Definition.} The map $h$ is called {\it transversal} at the point $x$
if the map $h_U$ is transversal to the stratum containing the point
$h_U(x)$. The map $h$ is called {\it transversal} if it is transversal at
every point $x\in M$.

Note that the transversality condition does not depend on the trivialization of the bundles 
chosen over $U$. Thom's transversality theorem implies that a generic
map $h$ is transversal at every point $x$.

\medskip
\noindent  
{\bf Corollary 1.} {\sl If a map $h$ is transversal then for any growth vector $r$ its  degeneracy locus
$\Omega_{r}$ is a closed (possibly empty) sub-variety of $M$. The dual $\Bbb Z_2$-cohomology
class given by the intersection index with the smooth part of
$\Omega_{r}$
is well-defined and given by the analog of formula (8) where the Chern classes substituted
by the corresponding Stiefel-Whitney classes.} 

{\em Proof.} One should follow  step-by-step the proof of Proposition 10.2
of \cite {F2}. In fact, one can show that the proof of formula (8) can be
reduced to the following basic results:

1) the axioms of the Chern classes the most important of which being the
Whitney formula $c(C\oplus F)=c(E)c(F)$.

2) the fact that the cohomology ring of $\Bbb CP^n$ is given by $H^*(\Bbb C P^n)=\Bbb
Z[c_1]/c_1^{n+1},$ where $c_1$ is the first Chern class of the
tautological bundle over $\Bbb CP^n$.

3) the construction of the Gysin map $\phi_*:H^*(X)\to H^*(Y)$ for the proper
map $\phi:X\to Y$ of smooth manifolds, and the Gysin formula
$\phi_*(\phi^*a\cup b)=a\cup \phi_*(b), a\in H^*(X), b\in H^*(Y).$

4) the relation  
$$p_*(1/(1-c_1(S))=c^{-1}(E^*),$$
where  $\pi:E\to M$ is any complex vector bundle, $p:P\to M$ is its
projectivization, and finally $S$ is the natural tautological linear sub-bundle in
the bundle $p^*\pi$ over $P$.

All these assertions have their real analogs with $\Bbb Z$-cohomology
 replaced by $\Bbb Z_2$-cohomology, the Chern classes by the
Stiefel-Whitney classes,
$\Bbb CP^n$ by $\Bbb RP^n$ etc. The proof of the corollary follows. \qed

\subsection {Application to sub-bundles}

If we drop the restrictions (7) then for a given n-sub-bundle $\CV\subset
TM_m$ and a given growth vector
$r=(r_1=n\le r_2\le ... \le r_k=m)$ the degeneracy locus
 $\Si_{r}$ is the subset
$\Omega_{r}\subset M$ for the map
$$\Phi_k: \fL^1(\CV)\subset...\subset \oplus_{i=1}^k\fL^i(\CV)\to TM.$$

Let us denote $L_j(\CV)=\oplus_{i=1}^j\fL^i(\CV)$ and
$\pa(n,j)=\dim(\oplus_{i=1}^j\fL^i(\CV))=\sum_{i=1}^jd(n,i).$
In order to apply Fulton's formula (8)
we must get rid of the redundant
subspaces, i.e. those subspaces whose rank conditions are automatically satisfied due to the rank conditions imposed on the previous subspaces. (In other words, the $i$th subspace is redundant if 
$\pa(n,i)-r_i=\pa(n,i-1)-r_{i-1}$.)  
We call by a {\it reduced index set} $I=(i_1,...,i_l)$
 the maximal subset of indices for which the corresponding $r_{i_j}$ satisfy
the conditions
$$r_{i_1}<...<r_{i_l}<m,\,
0<\pa(n,i_1)-r_{i_1}<\pa(n,i_2)-r_{i_2}<...<\pa(n,i_l)-r_{i_l},$$
i.e. both ranks and coranks are strictly increasing.

One gets the following Young diagram
$\la(r)=(p_1(r)^{m_1(r)},...,p_l(r)^{m_l(r)})$ where
$$p_1(r)=\pa(n,i_l)-r_{i_l},\,p_2(r)=\pa(n,i_{l-1})-r_{i_{l-1}},...,
\,p_l(r)=\pa(n,i_1)-r_{i_1},\eqno(9)$$
$$m_1(r)=m-r_{i_l},\,
m_2(r)=r_{i_l}-r_{i_{l-1}},...,\,m_l(r)=r_{i_2}-r_{i_1}.$$

Its dual diagram is
$\mu(r)=(q_1(r)^{n_1(r)},...,q_l(r)^{n_l(r)})$ where
$$q_1(r)=m-r_{i_1},...,q_l(r)=m-r_{i_l},$$
$$n_1(r)=\pa(n,i_1)-r_{i_1},...,n_l(r)=(\pa(n,i_l)-r_{i_l})-
(\pa(n,i_{l-1})-r_{i_{l-1}}).\eqno(10)$$

Finally, we set $cd({r})=\vert \la(r) \vert=\vert \mu(r) \vert$ is the
area of either of these Young diagrams.

Analogously, $s\la(r)=m-r_{i_1}$,
$s\mu(r)=\pa(n,i_l)-r_{i_l}$
and
$$\rho_r(i)=max\{s\in [1,l]: i\le
m-r_{i_s}=m_1(r)+...+m_{l+1-s}(r)\},\,i=1,...,s\la(r)$$
$$\rho_r^\prime(i)=min \{s\in [1,l]: i\le
\pa(n,i_s)-r_{i_1}=n_1(r)+...+n_{s}(r)\},\,i=1,...,s\mu(r) .$$

\noindent 
{\bf Definition.} The number $cd(r)$ is called the {\it expected
codimension} of $\Si_r$.
\bigskip\bigskip

\noindent 
{\bf Main result.}  {\sl If $\text{codim} (\Si_r(\CV))$ coincides with its
expected
 codimension $cd(r)$ then the $\Bbb Z_2$-cohomology class
$[\Si_r]_{\Bbb Z_2}$ of the base manifold 
 $M$ dual to $\Si_r$ is given by
$$[\Si_r]_{\Bbb
Z_2}=\det(w_{\la_{i}(r)-i+j}(L(\CV)_{\rho_r(i)}^*-TM^*))_{1\le i,j\le
s\la(r)}=$$ 
$$=\det(w_{\mu_{i}(r)-i+j}(TM-L(\CV)_{\rho_r^\prime(i)}))_{1\le i,j\le
s\mu(r)},\eqno(11)$$
where $w_l$ are the Stiefel-Whitney classes.}
\bigskip

\noindent 
{\bf Examples,} compare \cite {KMS}. Consider a generic 2-sub-bundle in
$TM_4$. There are 3 possible
non-maximal growth vectors $(2,2,4),\,(2,3,3),\,(2,2,3,4)$. The coincidence
of the actual and the expected codimensions in this case follows from the
normal forms in \cite {Zh}. (In the case
$(2,2,2,...)$ the codimension is $\ge5$.)

1) Case $r=(2,2,4)$. The reduced index set is $I=\{2\}$, i.e. we have to
consider
only the map $\Phi_2: L_2(\CV)\to TM$ of the usual bundles and determine the
locus of points where $rk(\Phi_2)\le2$. One has
$rk (L_2(\CV))=3,\,rk(TM)=4$, $\la(r)=(1^2),\,s\la(r)=2,\,cd(r)=2.$
Finally, $\mu(r)=2,\,s\mu(r)=1$ and $\rho^\prime(1)=1.$
Therefore,
$$[\Si_{(2,2,4)}]_{\Bbb Z_2}=w_2(TM-L_2(\CV))=w_2(M)+w_2(\CV)+w_1^2(\CV).$$

2) Case $r=(2,3,3)$.  The reduced index set is $I=\{3\}$, i.e. we have to
consider
only the map $\Phi_3: L_3(\CV)\to TM$ of the usual bundles and determine the
locus of points where $rk(\Phi_3)\le3$. One has
$rk (L_3(\CV))=5,\,rk(TM)=4$, $\la(r)=(2^1),\,s\la(r)=1,\,cd(r)=2.$
Finally, $\mu(r)=(1^2),\,s\mu(r)=2$ and
$\rho^\prime(1)=1,\,\rho^\prime(2)=1.$
Therefore,
$$ [\Si_{(2,3,3)}]_{\mathbb Z_2}  =
\left\vert \begin{array}{cc} w_1(TM-L_3(\CV))&w_2(TM-L_3(\CV))\\
w_0(TM-L_3(\CV))&w_1(TM-L_3(\CV))\end{array}\right\vert=$$
$$=w_1^2(M)+w_1^2(\CV)+w_2(M)+w_1(M)w_1(\CV).$$

3) Case $r=(2,2,3,4)$. The reduced index set is $I=\{2,3\}$, i.e. we have
to consider  the map $\Phi: L_2(\CV)\subset L_3(\CV)\to TM$.
One has
 $\la(r)=(2,1),\,s\la(r)=2,\,cd(r)=3.$  Now,
 $\mu(r)=(2,1),\,s\mu(r)=2$ and $\rho^\prime(1)=1,\,\rho^\prime(2)=2.$
Therefore,
$$[\Si_{(2,2,3,4)}]_{\mathbb Z_2}  =
\left\vert \begin{array}{cc} w_2(TM-L_2(\CV))&w_3(TM-L_2(\CV))\\
w_0(TM-L_3(\CV))&w_1(TM-L_3(\CV))\end{array}\right\vert=$$
$$=w_1(M)w_2(M)+w_2(M)w_1(\CV)+w_1^3(\CV)+w_3(M).$$

The above answers are obtained by the standard manipulations with the total
Chern
class of $TM$ and $\fL^k_n$.

\section {Transversality property for sub-bundles
and some general properties  on $f\CV$}\label{sec:4} 

In order to be able to apply formula $(11)$  to sub-bundles one needs to show 
that a certain transversality property is valid for the map 
$\Phi_k:\oplus_{i\leq k}\fL^i(\CV)\to TM$, see \S 2. This condition 
can be formulated as follows. The total space $Hom(\oplus_{i\leq 
k}\fL^i(\CV),TM)$ has a natural stratification according to different 
degenerations of the growth vector. The transversality property says 
that the section of the above bundle determined by the map $\Phi_k$ is 
transversal to each stratum of this natural stratification.  

   Naturally one wants to know if the transversality property is 
   valid for generic $n$-dimensional sub-bundles $\CV\subset TM$. The 
   conjecture stated below claims that this is indeed the case. 
   (Up to codimension $m-\sqrt{m}$ a similar statement is shown to 
   be valid in \cite {GV}.) 
    
   Since the transversality property is essentially local let us 
   formulate it in local terms.

\subsection {Local problem}

Take $M=\Bbb R^m$ with a fixed system coordinates $x_1,...,x_m$
  and consider the set $\Omega^0$ of germs of $n$-sub-bundles in $\Bbb R^m$ such
 that for any $\CV\in \OO$ the subspace $\CV(0)$ at the origin is spanned by
 $\frac{\pa}{\pa x_1},..,\frac{\pa}{\pa x_n}$.
The set $\OO$ can be identified with the set of
all $n$-tuples
of vector-fields $v_1,...,v_n$ of the form
$v_i=\frac{\pa}{\pa x_i}+\sum_{j=n+1}^ma_{i,j}(x_1,...
,x_m)\frac{\pa}{\pa x_j}$. Indeed, fixing
 the standard Eucledian structure on $\Bbb R^m$
we can  uniquely lift the vector fields $\frac{\pa}{\pa x_1},...,
\frac{\pa}{\pa x_n}$ to any sub-bundle
$\CV\in \OO$ and get the $n$-tuple of vector fields $v_1(\CV),...v_n(\CV)$
with the above properties.

\medskip 
\noindent 
{\bf  Remark.} For each $k$ we have the derived map
$\Psi_k: \;\CV\to Fl_k(\CV)$ where $Fl_k(\CV)=(\CV=\CV_1\subset
\CV_2\subset ... \subset \CV_k)$.  Each
$Fl_k(\CV)$
 is the image under the
 canonical map of the $\CA$-module $N_k(\CV)$, see  \S 2.
Fixing the standard Euclidean structure, the $n$-tuple of vector fields
$\pa_i$
 and the Hall basis we obtain the standard set of sections for all $\CV\in
\OO$ and in
all $N_k(\CV)$. This gives us a non-canonical isomorphism of $\CV_1$ and
$\CV_2$ and $N_k(\CV_1)$ and $N_k(\CV_2)$ as the $\CA$-modules for any sub-bundles
$\CV_1$, $\CV_2$ and any positive integer $k$. Localizing we can identify the jets
 of the sub-bundle $\CV$ with the jets of the $n$-tuples of vector fields $v_1(\CV),...,
v_n(\CV)$ and the jets of $Fl_k(\CV)$
with the jets of the $\pa(n,k)$-tuples of vector fields obtained from
$v_1(\CV),...,v_n(\CV)$ by applying the commutations prescribed by the
elements in the chosen Hall basis. (Recall that
$\pa(n,k)=\sum_{i=1}^kd(n,i)$.)
 The map $\Psi_k$ induces the
 well-defined map  $\Psi_k^i: j^{k+i}(\CV)\to j^i(Fl_k(\CV))$ of the corresponding jets.

\medskip 
\noindent 
{\bf Remark.} The $0$-jet of $Fl_k(\CV)$ can be represented by a
$m\times  \pa(n,i)$-matrix of the form
$$\left(\begin{array}{cccccccccccccc} &  \vdots & n&\vdots& m-n &\vdots & d(n,2)  &\vdots & d(n,3)
          &\vdots  & ... &\vdots  & d(n,k) &  \vdots\\
           \hdotsfor{14}\\
          n & \vdots & \one && 0   & &    0   & &    0   &  &
             ... & &    0   & \vdots  \\
         m-n&  \vdots & 0 && 0  &   &    * &  &    *  & & ... & & *  &
\vdots \\
         \hdotsfor{14}
\end{array} \right) \eqno(12)$$

Here $\one$ and $0$ denote the identity and the zero matrices of the sizes given in
the first row and the first column and $'*'$ stands for arbitrary real entries.

\medskip 
\noindent 
{\bf  Further notation.} Let  $Mat^0(n,m,k)$ denote the subset of all
$m\times\pa(n,k)$-matrices of the above form (12) and
 $J(n,m,k)$ denote the space of $k$-jets of $\CV\in \Omega^0$.
The space $J(n,m,k)$ is isomorphic to the space of all $n(m-n)$-tuples of
polynomials in $m$ variables of degree $\le k$, see above. 
 Obviously, $\dim J(n,m,k)=(m-n)n\binom {m+k}{k}$ and
 $\dim Mat^0(n,m,k)=(m-n)\sum_{i=2}^{k}d(n,i)$.
\medskip
The main object of the remaining part of the paper is the polynomial map of affine spaces
$$\Psi_k^0: J(n,m,k-1)\to Mat^0(n,m,k).$$
\medskip 
The space $Mat^0(n,m,k)$ has the following natural stratification. We 
start with the obvious  inclusions $Mat^0(n,m,1)\subset Mat^0(n,m,2)\subset ...
\subset Mat^0(n,m,k)$. Now fixing the growth vector $r=(r_1=n\le
r_2\le ... \le r_k\le \pa(n,k))$ we define the
 subset $Mat_r^0(n,m,k)
\subset Mat^0(n,m,k)$ of all matrices restriction of which to $Mat^0(n,m,i)$
have rank $\le r_i$ for all $i=1,..,k$.
Rather obviously, the codimension of $Mat_{r}^0(n,m,k)$ in $Mat^0(n,m,k)$
equals $cd(r)$.

Finally we are in position to formulate 
the required transversality property. 

\bigskip 
\noindent 
{\bf Main Conjecture.} {\sl The subset of points $x$ in $J(n,m,k)$ such that the map 
$\Psi_k^0$ is non-transversal to the stratum containing the point 
$\Psi_k^0(x)\in Mat^0(n,m,k)$ has the codimension in $J(n,m,k)$ 
strictly exceeding $m$. }

\bigskip 
 Thom's transversality theorem implies that the validity of the above 
 conjecture immediately leads to the transversality assumption for 
 generic $n$-dimensional sub-bundles in $TM$.
 
  We were unable to prove the above conjecture in its complete 
  generality but we were able to settle a number of cases given 
  below. 
   
\medskip
\noindent    
  {\bf Transversality Theorem.} {\sl The required transversality property holds either  for
$n\ge3$ and  $m\le \frac{n(n+1)(2n+1)}{6}=\pa (n,3)$ or for $n=2$ and
$m\le 8=\pa (2,4)$. Namely, for $k=2,3$ and any $m \ge n$ the map $\Psi^0_k$ is a
submersion. Also, for $n=2$ and any $m\ge n$ and $k=4$ the map $\Psi_4^0$ is a
submersion. } 
\medskip

\subsection {Violation of transversality for big codimensions}

We will finish this section by pointing out that the behavior of the image
$\Psi^0_k(J(n,m,k-1))$ w.r.t the natural rank stratification of the space
$Mat^0(n,m,k)$ is highly nontrivial. The next 2 statements show that
one can only hope that  transversality holds only
for the strata of  relatively small codimension as stated in the main conjecture. 

\medskip 
\noindent 
{\bf Lemma 1.} {\sl For any fixed $m\ge n$ and for
 $k>\text { const } log_{\binom n2}\binom {2m}{m}$) one has $\dim J(n,m,k-1)
< \dim Mat^0(n,m,k)$ and therefore $\Psi_k$ is not surjective.} 

{\em Proof.} This is simply the dimension count since $\dim J(n,m,k-1)=(m-n)n
\binom {m+k-1}{k-1}$ and  $\dim Mat^0(n,m,k)=(m-n)\sum_{i=2}^{k}d(n,i)$
where $kd(n,k)=\sum_{j\vert k}\mu(j)n^{k/j}$ with $\mu(j)$ denoting
the M\"obius function. \qed

\medskip 
\noindent 
{\bf Lemma 2.} {\sl For any $m\ge n\ge 3$ and $k\ge 4$ the map $\Psi_{k}$ is
 never onto. The same holds for $n=2$ and $k\ge 5$.}Ê

{\em Proof.}  Consider a matrix
in $Mat^0(n,m,4)$ of the form

$$\begin{array}{cccccccccccc} &  \vdots & n&\vdots&m-n &\vdots & d(n,2)  &\vdots & d(n,3)
          & \vdots  & d(n,4) &  \vdots\\
           \hdotsfor{12}\\
          n & \vdots & \one && 0   & &    0   & &    *   &
          &    \text{full}   & \vdots  \\
         m-n&  \vdots & 0   &   &   0 && 0 &  &    *  &  & \text{rank}  &
\vdots \\   \hdotsfor{12}
\end{array}.$$

Such a matrix is never in the image of $\Psi_4$ since the $4$-th homogeneous
component $H^4$ of the Hall basis  contains elements of the form
$((v_i,v_j),(v_p,v_q))$ which vanish since the first commutators of all basic
vector fields vanish. (For $n=2$ the same effect happens for
$k\ge 5$.)  \qed

\section {Transversality theorem and defect vectors}\label{sec:5}

This section is devoted to the proof of transversality theorem of the previous section as well as to the description of generic defect vectors in our situation. 
Namely, let 
$w=(((v_{i_1},v_{i_2})(v_{i_3},v_{i_4}))...((v_{i_{k-4}},v_{i_{k-3}})(v_{i_{
k-1}},v_{i_{k}})))
\in H^k$ be any element of length $k$ in the Hall basis.

\medskip 
\noindent 
{\bf Definition.}  By {\it the depth $dp_w(v_{i_j})$} of any variable
$v_{i_j}$ in $\,w$ we
denote the difference between the number of opening and closing parentheses
preceding
$v_{i_j}$. By the {\it depth $dp(w)$} of
$w$ we mean $\max_jdp_w(v_{i_j})$.

 Obviously, for any $j$ one has $dp_w(v_{i_j})\le dp(w)$
and there exist at least 2 different values of $j$ for which
$dp_w(v_{i_j})= dp(w)$.

Let $H^p$ denote the $p$th homogeneous component of $H$ and $H^{p,q}$ be its
subset containing all elements of depth $q$. Obviously, $q$ can vary between
$[\log_2(p)]+1$ and $p$. Let us restrict the map $\Psi_k^0$ to the union
$\bigcup_{j=1}^kH^{j,j}$.

\medskip 
\noindent 
{\bf Lemma 3.} {\sl For any $k$ and $m\ge n$ the restricted map
$$\Psi_k^{res}:J(n,m,k-1)\to Mat^0_{res}(n,m,k)$$
 is a submersion where $Mat^0_{res}(n,m,k)\subset Mat^0(n,m,k)$ contains only
rows corresponding to the elements from $\bigcup_{j=1}^k H^{j,j}$.}Ê

{\em Proof.} We use induction on $k$.

Base of induction, $k=2$. One has $(v_i,v_j)_l=a_{j,i,l}-a_{i,j,l}$ where
$(v_i,v_j)_l$ is the value of the $l$th component of the commutator
$(v_i,v_j)$ at the
origin, $a_{j,i,l}$ is the coefficient at $x_i$ in the $l$th component of
$v_j$ and
$a_{i,j,l}$ is the coefficient at $x_j$ in the $l$th component of $v_i$.
Since $a_{j,i,l}$
and $a_{i,j,l}$ are independent parameters in the pre-image the map
from $J(n,m,1)$ to $Mat^0_{res}(n,m,2)$ is, obviously, a submersion.

Step of induction. Assume that the statement is proved for $j<k$.
All elements in $H^{k,k}$ are in 1-1-correspondence with the
 commutators $(v_{i_1}(v_{i_2}$ $(...(v_{i_{k-1}},$ $v_{i_{k}}))...)$ where
$i_1\ge i_2\ge ...\ge i_{k-1}<i_{k}$, see def. of the Hall basis. (In
particular,
the number $\sharp_k(n)$ of such elements equals
$\sum_{j=0}^{n-1}\binom{j+k-2}{j}j$. For
example, $\sharp_2(n)=\frac{n(n-1)}{2},\quad $ $
\sharp_3(n)=\frac{n(n-1)(n+1)}{3},
\quad $ $\sharp_4(n)=\frac{(n-1)n(n+1)(n+2)}{8},\quad
$ $\sharp_5(n)=\frac{n(n-1)^2(n+2)(n+3)}{30},\quad
$ $\sharp_6(n)=\frac{n^2(n^2-1)
(n^2+n+2)}{16}$.) 

One has
$$(v_{i_1}(v_{i_2}(...(v_{i_{k-1}},v_{i_{k}}))...)_l=\text { const }
+\frac {a_{i_{k},i_1i_2...i_{k-1},l}}{c(i_1,...,i_{k-1})}-
\frac {a_{i_{k-1},i_1i_2...i_{k-2}i_{k},l}}{c(i_1,...,i_{k-2},i_k)}.$$
Here $(v_{i_1}(v_{i_2}(...(v_{i_{k-1}},v_{i_{k}}))...)_l$ is the value at
$0$ of
the $l$th component of the commutator vector field corresponding to the
 considered element of $H^{k,k}$.
The term  $\text{ 'const' }$ in the r.h.s. depends only on the $(k-2)$nd
jet of the basic
vector fields $v_1(\CV),..,v_n(\CV)$; $c(i_1,...,i_{k-1})$ and
$c(i_1,...,i_{k-2},i_k)$ are
factorial type constant factors. Finally, in the above formula we denote by
$a_{i_{k},i_1i_2...i_{k-1},l}$
(resp. $a_{i_{k-1},i_1i_2...i_{k-2}i_{k},l}$)  the variable coefficient at
the product
$x_{i_1}x_{i_2}\dots x_{i_{k-1}}$ in the $l$th component
of the basic vector field $v_{i_k}$ (resp. at the product
$x_{i_1}x_{i_2}\dots x_{i_{k-2}}x_{i_k}$
in the $l$th component of the basic vector field $v_{i_{k-1}}$).
One can easily check that  the variable coefficient
$a_{i_{k},i_1i_2...i_{k-1},l}$ appears
 only in the commutator
$(v_{i_1}(v_{i_2}(...(v_{i_{k-1}},v_{i_{k}}))...)_l$. 
 Since all these variable coefficients are independent parameters in the
pre-image we get
the  submersion of the space $\Omega_{\mathfrak J}$ on the space
$\widetilde{Mat}^0_{res}(n,m,k)$.
Here  $\Omega_{\mathfrak J}$ is the subspace of all sub-bundles in $\Omega$ with
some arbitrary fixed $(k-2)$-jet $\mathfrak J$
and  $\widetilde{Mat}^0_{res}(n,m,k)$ is the last block in
${Mat}^0_{res}(n,m,k)$. \qed

\medskip
\noindent 
{\bf Corollary 2.} {\sl The above Lemma implies the transversality theorem of the previous section.}Ê 

{\em Proof.}  Indeed, in  
both cases $n=2,\; m\le 8$ and $n\ge 3,\; m\le  \frac{n(n+1)(2n+1)}{6}=\pa (n,3)$ one has $H^k=H^{k,k}$.
\qed

\subsection {Description of generic defect vectors}

We give a  {\sl conjectural} description of 
all possible degenerations of the associated flag bundle $f\CV$ which might occur for a generic 
$n$-dimensional bundle $\CV$  in $TM$.

\medskip 
\noindent 
{\bf Definition.} For any given growth vector $r=(r_1=n,...,r_k=m)$
  we call the vector
 $def(r)=(\pa(n,1)-r_1=0; (\pa(n,2)-r_2)-(\pa(n,1)-r_1);...;
(\pa(n,k-1)-r_{k-1})-(\pa(n,k-2)-r_{k-2});0)$  {\it the defect vector
of $r$}.
(Up to the reduction of redundant indices this definition coincides with
$n_i$'s in $(9)$.)
A stratum
$St_r$ is called
 {\it admissible} {\it (resp. potentially admissible) }  if $\text {codim
}St_r\le m$
(resp. $cd(r)\le m$) and {\it bounding} {\it (resp. potentially
admissible)} if it is not
admissible (resp. potentially admissible) but there is no non-admissible
(resp. potentially nonadmissible) stratum containing
$St_r$ in its closure.

\medskip 
\noindent
{\bf Lemma 4.} {\sl For any $n\ge 3$ and any $k\ge 1$ one has
 $d(n,k+1)>\pa(n,k)=\sum_{j=1}^k d(n,j)$. For $n=2$ and any $k\ge 1$ one has
$d(2,k+1)+d(2,k+2)>\pa(2,k)=\sum_{j=1}^k d(2,j)$.} 

{\em Proof.} We will consider only the case $n\ge3$. We construct for each
element $a\in H^j,\;j\le k$ the unique element in $H^{k+1}$, i.e. embed
$\bigcup H^j_{j=1,...,k}$ into $H^{k+1}$, thus proving that
$d(n,k+1)>\pa(n,k)$. Recall that the Hall family is linearly ordered and
 each  element in the Hall family has the unique
representation in the form $(a(bc))$ where $a,b,c$ satisfy the conditions:
$a,b,c,bc$ belong to $H$, $a\ge b$ and $b<c$. For all $j<[\frac{k+1}{2}]$
we associate to any element $a\in H^j$ the element
$(a(v_1(v_1(...(v_1,v_2)...))))\in H^{k+1}$.
Now let $j>[\frac{k+1}{2}]$ and $a$ be some element in $H^j$. Then
$a=(f,g)$ where $f<g$.
Assume additionally that $2lng(f)+lng(g)\le k+1$. Then we associate to $a$
the element
$(h(f,g))$ where $h$ is the maximal element in $H^{k+1-lng(h)-lng(g)}$.
(Note that, by
definition, $h\ge f$.)
If $2lng(f)+lng(g)> k+1$ then we choose for each $a=(f,g)\in H^j$ the element
$(f(h,g))$ where $h$ is the minimal element in $H^{k+1-lng(h)-lng(g)}$. The
last case
to consider is $H^l$ when $k=2l-1$. One has that $\Lambda^2(H^l)$ is embedded
into $H^{2k}$ and under the assumption $n\ge 3$ one has $\dim \Lambda^2(H^l)\ge
\dim H^l$. Combining all choices together we obtain a set-theoretic
embedding of
$\bigcup_{j=1}^k H^j$ into $H^{k+1}$.  The
result follows. More detailed consideration shows that
$d(2,k+1)+d(2,k+2)>\pa(2,k)$.
\qed

\medskip 
\noindent 
{\bf Remark.} The transversality property is equivalent to 
showing that the codimension of each potentially admissible (resp.
potentially bounding) stratum $St_r$ equals $cd(r)$ and therefore $St_r$
is, in fact, admissible (resp. bounding).

 Let us  enumerate the defect vectors of all potentially admissible
strata.

\medskip 
\noindent 
{\bf Lemma 5.} {\sl For $n\ge3$ the defect vectors of all potentially
admissible strata
are as follows. We assume that $\pa(n,p)\le m <\pa(n,p+1)$.} 

$$  \begin{array}{ccccccccccccc}
      \text {case/posit} &\vdots     &1      &2     &... &l-1      &l      &l+1
&... &p-1   &p        &p+1 & \vdots\\
\hdotsfor{13}\\
       a)&\vdots&         0      &0     &...  & 0      &1       &0       &...
&0 &\chi    &0     &\vdots \\
       b)&\vdots&         0      &0     &...   &0      &1       &0       &...
&  0   &0  & \chi& \vdots\\
       c)&\vdots&         0      &0     &...   &0      &0       &0       &...
&0&\chi   &\nu  & \vdots \\
\hdotsfor{13}
\end{array}$$
where $l<p$ and
 the following restrictions are satisfied for each of the above cases

$$\begin{cases}
a)\quad (m-\pa(n,p)+1+\chi)\chi\le
\pa_{l}-1,\,\chi\ge0; \\
b) \quad (m-\pa(n,p+1)+1+\chi)\chi\le
\pa_{l}-1,\,\chi+1+m-\pa(n,p+1)\ge0; \\
c) \quad (m-\pa(n,p)+\chi)\chi+(m-\pa(n,p+1)+\chi+\nu)\nu\le m,\\
(m-\pa(n,p+1)+\chi+\nu)\ge0.
\end{cases}
$$

{\em Proof.\;}  If the rank of the image of some
$L_l(\CV)=\oplus_{i=1}^l \mathfrak{Lie}^i(\CV)$ where $l<p$ drops then it can drop
exactly by 1. Indeed, assume that it drops by at least 2 then using $(10)$
we get for the dual diagram $\mu_{r}$ that $q_1\ge m-\pa(n,l)+2$ and
$n_1\ge 2$.
But by lemma 5.4. one has $2(m-\pa(n,l)+2)>m$ which contradicts to the
assumption
$cd(r)\le m$. Analogously, if $L_l(\CV),\,l<p$  has the corank equal to 1
then
the corank can possibly increase again (becomes $>1$) only for 
either $L_p(\CV)$ or
$L_{p+1}(\CV)$, see more detailed description below. Indeed, assume that the
corank  drops by 1 at $L_{l_1}(\CV)$ and by 2 at $L_{l_2}(\CV)$ where
$l_1<l_2<k$. Then by $(10)$ one has $q_1=m-\pa(n,l_1)+1,\,n_1=1,\,
q_2=m-\pa(n,l_2)+2,\,n_2=1$.
But again by 5.4. one gets $q_1+q_2>m$ which contradicts to $cd(r)\le m$.
Therefore
the second rank drop can only occur either at position $p$ or $p+1$. If we
have that
$L_l(\CV),\,l<p$ has corank 1 then further corank drops at both positions
$p$ and
$p+1$ are simultaneously impossible, i.e. only 1 extra drop is allowed. 
Indeed, assume that we have
coranks $1,2,3$ at positions $l,p,p+1$ resp. Then by $(10)$ one gets
$q_1=m-\pa(n,l)+1,n_1=1$; $q_2=m-\pa(n,p)+2,n_2=1$; $q_3=1,
n_3=\pa(n,p+1)-m-1$.
But $q_2n_2+q_3n_3=m-\pa(n,p)+2+\pa(n,p+1)-m-1=d(n,p)+1$. Again, by 5.4.
one gets $q_1n_1+q_2n_2+q_3n_3>m$. Thus we are left with the cases a) and b)
if the rank drops in some position prior to $p$. The above list of inequalities
 follows from the expressions for the terms in the dual Young diagram
$\mu=(q_1^{n_1},q_2^{n_2})$ given below. In the case a) one gets
$q_1=m-\pa(n,l)+1, n_1=1, q_2=m-(\pa(n,p)-1-\chi),n_2=\chi$. In the case b)
one gets $q_1=m-\pa(n,l)+1, n_1=1,
q_2=m-(\pa(n,p+1)-1-\chi),n_2=\chi$. Finally, in the case c) one gets
 $q_1=m-\pa(n,p)+\chi,n_1=\chi, q_2=m-\pa(n,p+1)+\chi+\nu,n_2=\nu$.
The inequalities express
 the condition $cd(r)\le m$ and the condition that the second rank drop
actually
occurs. One can easily show
that a) are b) are pairwise excluding, i.e. for a given pair $(n,m)$ either
inequalities for  a) or
or b) can be satisfied under the assumption that $\chi>0$.
\qed

\medskip 
\noindent 
{\bf  Corollary 3. } {\sl The defect vectors for all potentially bounding strata
have one of the following forms

$$ \begin{array}{ccccccccccccccccc} 
          & \vdots     &1      &2     &... &l_1-1    &l_1    &l_1+1 &...&l_2-1&l_2    &l_2+1 &... &p-1   &p        &p+1 & \vdots\\
       \hdotsfor{16}\\
       a)&\vdots&         0      &0     &...  & 0      &1       &0&... &...   &...    &...&... &0 &\chi    &0     &\vdots \\
       b)&\vdots&         0      &0     &...   &0      &1       &0&... &...   &...    &...&...&  0   &0  & \chi& \vdots\\
       c)&\vdots&         0      &0     &...   &...    &...     &...&... &...   &...    &...&...&0&\chi   &\nu  & \vdots \\
       d)& \vdots  &0       &0     &... &0       &1      &0      &...&... &1      &0     &...&...&0  &0       &\vdots  \\
       e)& \vdots &0        &0     &... &0       &2     &0       &...&... &...    &...  &...&...&0  &0       &\vdots\\
     \hdotsfor{16}
\end{array} $$

In the cases a)-c) (whose form coincides with that in Lemma~5) we additionally
require
that the inequalities from the formulation of Lemma~5 are violated for the considered values of
$\nu$ and $\chi$
but satisfied for any smaller nonnegative values of these variables.} 

{\em Proof.\;} Obtained by a simple case-by-case consideration.
\qed

\section {Final remarks} \label{sec:6}

\noindent 
{\bf 1.} The basic problem related to  the above topological obstructions is to what 
extent vanishing of these obstructions guarantees the existence of a maximal growth sub-bundle. 
The result of T.~Vogel on the existence of Engel structures on parallellizable $4$-manifolds brings a certain amount of optimism about this problem. For example, 

\medskip 
\noindent 
{\bf Problem.} Does every closed parallelizable $m$-dimensional manifold admits  a maximal growth distribution of rank $1<n<m$?  

\medskip 
\noindent 
{\bf 2.} The developed theory is  incomplete due to the
lack of proof of the transversality conjecture. Essentially, to
accomplish the proof one should only consider the defect
vectors listed in the end of \S 5. The authors are convinced that  progress in this
direction will  be intimately related to the detailed study of
the combinatorics of different Hall bases.

\medskip 
\noindent
{\bf 3.} The natural representation of $GL_n(\bC)$ in $\fL^i_n$
was extensively studied in the series of papers \cite {Thr}, \cite {We},
\cite {Kl} and some others. There exists an interesting formula
for the multiplicity of each irreducible representation (corresponding
to some Young diagram $\mu$) in $\fL^i_n$ analogous to the one for the
dimension
of $\fL^i_n$. This suggests the existence of a much more sophisticated theory
 of
characteristic classes for sub-bundles in the tangent bundle since there exist
many different natural filtrations in $\fL_n$ the most complicated of which
coming from
 the direct sum of irreducible $GL_n(\bC)$-representations. The 1st step in
this direction will be
to find analogs of Theorem 2.1 for these other filtrations and the second step
to find the analogous transversality theorems. As an example of
characteristic classes different from the ones considered above one
can try to calculate the characteristic classes related to the depth
filtration of $\fL_n$ introduced in \S 5. For this case there exists
a character formula in terms of plethysms analogous to $(5)$ suggested
to the authors by C.~Reutenauer but it is unfortunately rather
complicated.

\medskip 
\noindent 
{\bf 4.} Another natural question related to the transversality theorem is
to understand  for
strata of what corank in the space  $Mat^0(n,m,k)$ it holds and therefore to extend
the determinantal formulas to the case of non-generic sub-bundles or families
of sub-bundles.

\medskip 
\noindent 
{\bf Problem.} Generalize the transversality theory to the case of
non-generic sub-bundles.

\section {Appendix}\label{sec:7}

{\def\back{\catcode`\^=7\catcode`\_=8}
\def\ppar{\par\hskip0pt\hbox to 0pt{}}
\catcode`\^^M=13\let^^M=\ppar%
\catcode`\@=12\catcode`\#=12\catcode`\&=12
\catcode`\_=12\catcode`\^=12\catcode`\ =13
\parindent0pt\tt
(* MATHEMATICA program for calculation of Chern classes of
   homogeneous components {\back$\fL^k_n$} of free Lie algebra bundles
   up to the order                                            *)
order=4;

(* This function computes Chern character via total Chern class       *)
classtochar[cc_]:=
   (resc[n-t D[Log[cc],t],-1]/.t^k_->t^k/k!)+O[t]^(order+1)//ExpandAll;

(* This function calculates total Chern class via Chern character     *)
chartoclass[ch_]:=Exp[-Integrate[PolynomialQuotient[resc[ch,-1]/.
   t^k_->t^k k!,t,t],t]+O[t]^(order+1)]//ExpandAll;

(* This function makes rescaling     {\back$\eta\longmapsto(\eta)_l$}
*)
resc[eta_,l_]:=Normal[eta]/.t->l t;

(* Total Chern class of original bundle                               *)
class=1+Sum[c[i]t^i,{i,order}]+O[t]^(order+1);

(* Chern character of original bundle                                 *)
char=classtochar[class];

(* Calculates Chern character of    {\back$\fL^k$}    using Theorem 2.1.
*)
charofL[k_]:=
   (Plus@@((MoebiusMu[#]resc[char^(k/#),#]/k)&/@(Divisors[k])))//ExpandAll;

(* Total Chern class  of  {\back$\fL^k$}  as a series in $t$
*)
classofL[k_]:=chartoclass[charofL[k]];

}


\bigskip

\begin{thebibliography}{30}

\bibitem {At} M.~F.~Atiyah, {\em Vector fields on manifolds}, Arbeitsgemeinschaft f\"ur Forschung des Landes Nordheim-Westfalen, Heft 200, Westdeutscher Verlag, Cologne, 1970, 26 pp. 

\bibitem {Bou}  N.~Bourbaki, {Lie groups and Lie algebras}, Chapters 1-3, Springer-Verlag.

\bibitem {Ge1} V.~Ya.~Gershkovich, {\em On simplest Engel structures on 
$4$-manifolds}, Dynam. systems \& Appl., vol. 4,(1995),  274--294.

\bibitem {Ge2} V.~Ya.~Gershkovich, {\em On normal forms of distribution 
jets}, Lect. Notes in Math., vol. 1346, (1988), 77--98. 

\bibitem {Ge3} V.~Ya.~Gershkovich, {\em Exotic Engel structures on 
$\bR^4$}, Russian J. Math. Phys., vol. 3(2), (1995), 207--226. 


\bibitem {GV} V.~Ya.~Gerkshkovich and A.~M.~Vershik, {\em Nonholomorphic 
manifolds and nilpotent analysis},  J. of Geometry and Physics, vol.  
5(3), (1988), 407--452. 

\bibitem {F1} W.~Fulton, {Intersection theory},  Springer-Verlag, (1984).

\bibitem {F2} W.~Fulton, {\em Flags, Schubert polynomials, degeneracy
loci, and determinantal formulas}, Duke Math. J., 1992, vol. 65, 
(1992), 381--420.

\bibitem {FP} W.~Fulton and P.~Pragacz, {Schubert varieties and degeneracy loci}, Lecture Notes in Mathematics, vol. 1689, Springer-Verlag, Berlin, 1998, xii+148 pp. 

\bibitem{HH} F.~Hirzebruch, H.~Hopf, {\em Felder von Fl\"achenelementen in $4$-dimensionalen Mannigfaltigkeiten}, Math. Ann., vol. 136, (1958), 156--172. 

\bibitem {JM}ÊH.~Jacobowitz and G.~Mendoza, {\em Sub-bundles of the complexified tangent bundle}, Trans. Amer. Math. Soc., vol. 355(10), (2004), 4201--4222.   

\bibitem {Kl} A.~A.~Klyachko, {\em Lie elements in the tensor 
algebra},  Siberian Math. J., vol. 15, (1974), 914--920.

\bibitem {KMS}  M.~\`E.~Kazarian, R.~Montgomery, and B.~Z.~Shapiro, 
{\em Characteristic classes for the degenerations of  two-plane
fields in four dimensions}, Pacific J.Math., vol. 179(2), (1997), 
397--417. 

\bibitem {La}ÊD.~Laksov, {\em Remarks on Giovanni Zeno Giambelli's work and life. Algebra and geometry (1860--1940): the Italian contribution (Cortona, 1992).}, Rend. Circ. Mat. Palermo (2) Suppl. No. 36 (1994), 207--218.

\bibitem {Ma}ÊY.~Matsushita, {\em Fields of $2$-Planes on Compact Simply-Connected Smooth 
$4$-Manifolds}, Math. Ann., vol. 280, (1988), 687-689. 

\bibitem{Mont} R.~Montgomery, {A tour of subriemannian geometries, their geodesics and applications.} 
Math. Surveys and Monographs, vol. 91, AMS, Providence, RI, 2002, xx+259 pp. 

\bibitem {Re} C.Reutenauer,  {Free Lie algebras},  London Math. Soc., New 
Ser., vol. 7, (1993). 

\bibitem {Se} J-P.~Serre, {Lie algebras and Lie groups}, Benjamin, New 
York-Amsterdam, (1965). 

\bibitem {Th1} E.~Thomas, {\em Fields of Tangent $k$-Planes on Manifolds}, 
Invent. Math., vol. 3, (1967), 334--347. 

\bibitem {Th2} E.~Thomas, {\em Vector fields manifolds}, Bull. Amer. Math. Soc.,,  vol. 75, (1969), 643--683. 

\bibitem {Th3} E.~Thomas, {\em Fields of tangent $2$-planes on even-dimensional manifolds}, Ann. of Math., (2), vol. 86, (1967), 349--361. 

\bibitem {Thr} R.~M.~Thrall, {\em On symmetrized Kronecker powers and the 
structure of the free Lie ring},  Amer. J. Math., vol. 64,  (1942),  371--388. 

\bibitem {Vo} Th.~Vogel, {\em Existence of Engel structures}, math.GT/0411217, 57 pp. 

\bibitem {Vo1} Th.~Vogel, {\em A construction of Engel structures}, C.R.Math.Acad. Sci. Paris., 
vol. 340(1), (2005), 43--48.  

\bibitem {We}  F.~Wever, {\em \"Uber Invarianten in Lieschen Ringen}, 
 Math. Ann., vol. 120, (1949), 563--580. 

\bibitem {Zh} M.~Ya.~ Zhitomirskii, {\em Normal forms of germs of 2-dimensional 
distributions in $\R^4$}, Funct. Anal. Appl., vol. 24, (1990), 81--82. 


\end{thebibliography}
\end{document}